\documentclass[12pt]{article}

\usepackage{amssymb}

\usepackage{amsmath}

\usepackage[mathcal]{euscript}

\usepackage{amsthm}

\newcommand{\e}{\varepsilon}
\newcommand{\EE}{\mathsf{E}}
\newcommand{\PP}{\mathsf{P}}
\newcommand{\XX}{\mathbb{X}}
\newcommand{\YY}{\mathbb{Y}}
\newcommand{\DD}{\mathbb{D}}
\newcommand{\GG}{\mathbb{G}}
\newcommand{\NN}{\mathbb{N}}
\newcommand{\HH}{\mathbb{H}}

\begin{document}

\begin{center}
{\bf \large Reward Algorithms \vspace{1mm} \\ for Semi-Markov Processes}
\end{center} 

\begin{center}
{\large Dmitrii Silvestrov\footnote{Department of Mathematics, Stockholm University, SE-106 81 Stockholm, Sweden. \\ 
Email address: silvestrov@math.su.se} 
and Raimondo Manca\footnote{Department of Methods and Models for Economics, Territory and Finance,
University ``La Sapienza",  00161 Rome, Italy, E-mail address: raimondo.manca@uniroma1.it}}
\end{center}
\vspace{2mm}

Abstract:
New algorithms for computing power moments of hitting times and accumulated rewards of hitting type  for semi-Markov processes. The algorithms are based on  special techniques of sequential phase space reduction and recurrence relations connecting moments of rewards. Applications are discussed as well as possible generalizations of presented results and examples. \\
 


{\bf 1. Introduction}  \\

In this paper, we study recurrent relations for power moments of  hitting
times and accumulated rewards of hitting type for semi-Markov processes and present effective algorithms for computing
these moments. These algorithms are based on procedures of sequential of 
phase space reduction for semi-Markov processes.

Hitting times are often interpreted as transition times for different stochastic systems describing by Markov-type processes, for example,
occupation times or waiting times in queuing systems, life times in reliability
models, extinction times in population dynamic models, etc.  We refer to works by Korolyuk, Brodi and Turbin (1974), Kovalenko (1975), Korolyuk and Turbin  (1976, 1978), Courtois (1977), Silvestrov (1980b), Anisimov,  Zakusilo and Donchenko (1987), Ciardo,  Raymonf,  Sericola and Trivedi (1990), Kovalenko,  Kuznetsov and  Pegg (1997), Korolyuk, V.S. and Korolyuk, V.V. (1999),     Limni\-os and Opri\c{s}an  (2001, 2003), Barbu,   Boussemart and  Limnios (2004), Yin and Zhang (2005,  2013), Janssen and  Manca (2006,  2007),  Anisimov (2008),  Gyllenberg and Silvestrov (2008), D'Amico, Petroni and  Prattico (2013), and  Papadopoulou (2013).

In financial and insurance applications, the hitting times for semi-Markov processes can be also interpreted as  rewards  accumulated up to
some hitting terminating  time for a  financial or insurance contract. We refer here to works by D'Amico,  Janssen and Manca (2005),  Janssen and  Manca (2006, 2007), Stenberg, Manca and  Silvestrov (2006, 2007), Biffi,  D'Amigo,   Di Biase,   Janssen,  Manca and  Silvestrov (2008), Silvestrov, Silvestrova and Manca (2008), 
D'Amico and Petroni (2012),  Papadopoulou, Tsaklidis, McClean and  Garg (2012), D'Amico, Guillen and  Manca (2013),  and 
D'Amico, Petroni and  Prattico (2015).

Moments of hitting times also  play an important role in limit and ergodic
theorems for Markov type processes. As a rule, the first and second order moments  are used in conditions of theorems, higher order
moments in rates of convergence and  asymptotical expansions.  We refer here to works by Silvestrov (1974, 1980b, 1994, 1996), Korolyuk and Turbin  (1976, 1978),  Korolyuk, V.~S. and Korolyuk, V.~V. (1999), Koroliuk and Limnios (2005),  Anisimov (2008), Gyllenberg and Silvestrov (2008),  Hunter  (2005), Yin and Zhang (2005,  2013), Silvestrov and Drozdenko (2006) and Silvestrov, D. and Silvestrov, S. (2015). 

Recurrent relations, which link power moments of  hitting times for Markov chains have been first obtained for
Markov chains by Chung (1954, 1960). Further development have been
achieved by Lamperty (1963),  Kemeny and Snell (1961a, 1961b),  and  Pitman (1974a, 1974b, 1977), Silvestrov
(1980a, 1980b). Similar relations as well as description of these moments as minimal solutions
of some algebraic or integral equations were considered for  Markov chains and semi-Markov processes with discrete and
arbitrary phase spaces by Cogburg (1975), Nummelin (1984), and Silvestrov
(1980b, 1983a, 1983b, 1996) and Silvestrov,  Manca and  Silvestrova (2014). Analogous results  for mixed power exponential moments of first hitting
times for semi-Markov processes have been obtained  in Gyllenberg and Silvestrov (2008).

The paper includes five sections. In Section 2, we introduce Markov renewal processes, semi-Markov processes and define 
hitting times and accumulated rewards of hitting type. We also present basic stochastic relations  and recurrent systems of linear equations for power moments  of  these random functionals.  In Section 3, we describe a procedure of phase space reduction for semi-Markov processes and formulas for computing transition characteristics for reduced semi-Markov processes, We also prove invariance of hitting times and their moments with respect to the above procedure of phase space reduction. In Section 4, we describe a procedure of sequential phase space reduction for semi-Markov process and derive recurrent formulas for computing power moments of hitting times for semi-Markov processes. In Section 5, we present useful generalizations of the above results to real-valued  and vector accumulated rewards of hitting type, general hitting times  with hitting state indicators, place-dependent and time-dependents 
hitting times  and accumulated rewards of hitting type and give a numerical example  
for the corresponding recurrent algorithms for computing power moments of hitting times and accumulated rewards of hitting type for 
semi-Markov processes. \\

{\bf 2. Semi-Markov processes and  hitting times} \\

In this section, we introduce Markov renewal processes and semi-Markov processes.  We define also hitting times and accumulated rewards of hitting times, 
and give basic recurrent system of linear equations for  their  power moments, which are the main objects of our studies. \vspace{1mm}

{\bf 2.1. Markov renewal processes and semi-Markov processes}. Let  ${\mathbb X} = \{0, \ldots, m \}$ and $(J_n, X_n), n = 0, 1, \ldots$ be a Markov renewal process, i.e., a homogeneous Markov chain with the phase space ${\mathbb X} \times [0, \infty)$, an initial distribution $\bar{p} = \langle p_i = \PP \{J_0 = i, X_0 = 0 \} = \PP \{J_0 = i \}, i \in {\mathbb X} \rangle$ and transition probabilities,
\begin{equation}\label{semika}
Q_{ij}(t) = \PP \{ J_{1} = j, X_{1} \leq t / J_{0} = i, X_{0} = s \}, \ (i, s), (j, t) \in  {\mathbb X} \times [0, \infty).
\end{equation}

In this case,  the random sequence $\eta_n$ is also a homogeneous (embedded) Markov chain with the  phase space $\XX$ and the transition probabilities,
\begin{equation}\label{embed}
p_{ij} = \PP \{ J_{1} = j / J_{0} = i \} = Q_{ij}(\infty), \ i, j \in \XX.
\end{equation}

As far as random variable $X_n$ is concerned, it can be interpreted as sojourn time in state $J_n$, for $n = 1, 2, \ldots$.

We assume that the following communication conditions hold:
\begin{itemize}
\item[${\bf  A}$:] $\XX$ is a communicative class of states for the embedded Markov \mbox{chain $J_n$.}
\end{itemize}

We also assume that the following condition excluding instant transitions holds:
\begin{itemize}
\item[${\bf B}$:]  $Q_{ij}(0) = 0, i, j \in \mathbb{X}$.
\end{itemize}

Let us now introduce  a semi-Markov process,
\begin{equation}\label{sepras}
J(t) = J_{N(t)},  \ t \geq 0,
\end{equation}
where $N(t) = \max(n \geq 0: T_n \leq t)$
is a number of jumps in the  time interval  $[0, t]$, for $t \geq 0$, and
$T_n = X_1 + \cdots + X_n, \ n = 0, 1, \ldots$,
are sequential moments of jumps, for the semi-Markov process $J(t)$.

This process has the phase space $\XX$, the initial distribution $\bar{p} = \langle p_i = \PP \{J(0) = i \}, i \in {\mathbb X} \rangle$ and transition probabilities $Q_{ij}(t), t \geq 0,  i, j \in \XX$. \vspace{1mm}

{\bf 2.2. Hitting times and accumulated rewards of hitting type}. Let us also introduce moments of sojourn times,
\begin{equation}\label{expev}
e_{ij}^{(r)} = \EE_i   X_{1}^r  I( J_{1} = j) = \int_0^\infty t^r Q^{(\e)}_{ij}(dt), \ r = 0, 1, \ldots, i, j \in  \XX.
\end{equation}

Here and henceforth,  notations $\PP_i $ and $ \EE_i $ are used for conditional probabilities and expectations under
condition $J(0) = i$.

Note that,
\begin{equation}\label{expeva}
e_{ij}^{(0)} = p_{ij}, \ i, j \in  \XX.
\end{equation}

We assume that the following condition holds, for some integer $d \geq 1$:
\begin{itemize}
\item[${\bf C_d}$:]  $e^{(d)}_{ij} < \infty, \,  i, j \in \mathbb{X}$.
\end{itemize}

The first hitting time to state $0$ for the semi-Markov process $J(t)$ can be defined as,
\begin{equation}\label{hitta}
W_0 = \inf(t \geq  X_1: J(t) = 0) = \sum_{n = 1}^{U_0} X_n,
\end{equation}
where $U_0 =  \min(n \geq 1: J_n = 0)$ is the first hitting time to state $0$ for the Markov chain $J_n$.

The random variable $W_0$ can also be interpreted as a reward accumulated  on trajectories of Markov chain $J_n$ up to its first hitting to  
state $0$.  

The main object of our studies are power moments for the first hitting times,
\begin{equation}\label{mhitta}
E_{i0}^{(r)} = \EE_i W^r_0,   \  r  = 1, \ldots, d, \,   i \in \XX.
\end{equation}

Note that,
\begin{equation}\label{mhittama}
E_{i0}^{(0)} = 1, \,  i \in {\mathbb X}.
\end{equation}

As well known, conditions ${\bf A}$,  ${\bf B}$ and ${\bf C_d}$ imply that,
\begin{equation}\label{mhittav}
E_{i0}^{(r)}  < \infty, \  r =1, \ldots, d,  \ i \in \XX.
\end{equation}

In what follows,  symbol $Y \overset{d}{=} Z$ is used to denote that random variables or vectors $Y$ and $Z$  have the same distribution.

The Markov property of the Markov renewal process $(J_n, X_n)$ implies that following system of stochastic equalities takes
place for  hitting times,
\begin{equation}\label{inadopion}
\left\{
\begin{array}{ll}
W_{i, 0} \overset{d}{=}  X_{i, 1} I(J_{i, 1} = 0) + \sum_{j \neq 0}(X_{i, 1} + W_{j, 0})I(J_{i, 1} = j),  \\
 i \in \mathbb{X},
\end{array}
\right.
\end{equation}
where:  (a) $W_{i, 0}$ is a random variable which has distribution $\mathsf{P} \{ W_{i, 0} \leq t \} = \mathsf{P}_i \{ W_0 \leq t \}, t \geq 0$, for every $i \in \XX$;  (b) $(J_{i, 1}, X_{i, 1})$ is a random vector, which  takes values in space ${\mathbb X} \times [0, \infty)$ and has the distribution $P \{J_{i, 1} = j, X_{i, 1} \leq t \} = Q_{ij}(t),  j \in \XX, t \geq 0$, for every $i \in \XX$; (c) the random variables $W_{i, 0}$ and the random vector $(J_{i, 1}, X_{i, 1})$ are independent, for every $i \in \XX$.

By taking expectations in stochastic relations (\ref{inadopion}) we get the following system of linear equations
for expectations of hitting times  $E^{(1)}_{i0}, i \in \XX$, 
\begin{equation}\label{inadopiona}
\left\{
\begin{array}{ll}
E^{(1)}_{i0} = e^{(1)}_{i0}  + \sum_{j \in \XX, j \neq 0} e^{(1)}_{ij} + \sum_{j \in \XX, j \neq 0} p _{ij} E^{(1)}_{j0},  \\
i \in {\mathbb X}.
\end{array}
\right.
\end{equation}

In general,  taking  moments of the  order $r$  in stochastic relations (\ref{inadopion}) we get the following system of linear equations  for moments  $E^{(r)}_{i0}, i \in \XX$, for $r = 1, \ldots, d$,
\begin{equation}\label{inadopionba}
\left\{
\begin{array}{ll}
E^{(r)}_{i0} = f^{(r)}_{i0}   + \sum_{j \in \XX, j \neq 0} p_{ij} E^{(r)}_{j0},   \\
i \in {\mathbb X},
\end{array}
\right.
\end{equation}
where
\begin{equation}\label{inadop}
f^{(r)}_{i0}  = e^{(r)}_{i0} + \sum_{j \in \XX, j \neq 0} \sum_{l=0}^{r-1} \binom{r}{l}e^{(r-l)}_{ij} E^{(l)}_{j0}, \ i \in \XX.
 \end{equation}

The  system of linear equation  given in (\ref{inadopionba})  has, for $r = 1, \ldots, d$,  the same matrix of coefficients ${\bf I} - {\bf P}_0$, where ${\bf I}  = \| I(i = j) \|$ is the unit matrix and  matrix ${\bf P}_0 = \| p_{ij}I(j \neq 0) \|$.  

This is readily seen that matrix ${\bf P}^n_0 = \| \PP_i \{ U_0 > n, J_n = j \} \|$. Condition ${\bf A}$ implies that $\PP_i \{ U_0 > n, J_n = j \}  \to 0$ as $ n \to \infty$, for $i, j \in \XX$ and, thus, $\det({\bf I}  - {\bf P}_0) \neq 0$.

Therefore,  moments $E^{(r)}_{i0}, i \in \XX$ are the unique solution for the system of linear equation (\ref{inadopionba}), for
every $r = 1, \ldots, d$.

These systems have a recurrent character, since, for every $r 0 1, \ldots, d$,  the free terms $f^{(r)}_{i0} = f^{(r)}_{i0}(E^{(k)}_{j0}, j \neq 0, k = 1, \ldots, r -1), i \in \XX$ of the system  (\ref{inadopionba})  for moments 
$E^{(r)}_{i0}, i \in \XX$  are functions of moments $E^{(k)}_{j0}, j \neq 0, k = 1, \ldots, r -1$.

Thus, the systems  given in (\ref{inadopionba})  should be solved recurrently, for $r = 1, \ldots, d$.

This is useful to note that the above remarks imply that condition ${\bf  A}$ can be replaced by simpler hitting condition:
\begin{itemize}
\item[${\bf  A_0}$:] $\PP_i \{ U_0  < \infty \} = 1, \, i \in \XX$.
\end{itemize}

Let denote matrix $[{\bf I}  - {\bf P}_0]^{-1} = \| g_{i0j} \|$. The elements of this matrix have the following probabilistic
sense, $g_{i0j}  = \EE_i \sum_{n = 1}^{U_0} I(J_{n -1} = j), \, i, j \in \XX$.

The recurrent formulas for moments $E^{(r)}_{i0}, i \in \XX$ have the following form, for $r = 1, \ldots, d$,
\begin{equation}\label{inadopio}
E^{(r)}_{i0} = \sum_{j \in \XX} g_{i0j} f^{(r)}_{j0}(E^{(k)}_{l0}, l \neq 0, k = 1, \ldots, r -1), \, i \in \XX.
 \end{equation}

 This method of computing moments $E^{(r)}_{i0}, i \in \XX$ requires to compute the inverse matrix
 $[{\bf I}  - {\bf P}_0]^{-1}$.

In this paper, we propose an alternative method, which can be considered as a stochastic analogue of Gauss
 elimination method for solving of the recurrent systems of linear equations  (\ref{inadopionba}). \\

{\bf 3. Semi-Markov processes with reduced phase spaces} \\

 In this section, we describe an one-step algorithm for reduction of phase space for semi-Markov processes. We also give recurrent systems of linear equations for power moments of hitting times for reduced semi-Markov processes. \vspace{1mm}

{\bf 3.1. Reduced semi-Markov processes}.  Let us choose some state $k \in \XX$ and consider the reduced phase 
space $_k\XX = \XX \setminus \{ k \}$,  with the state $k$ excluded from the phase space $\XX$.

Let us define  the sequential moments of hitting the reduced space
$_k\XX$ by the embedded Markov chain $J_n$,
\begin{equation}\label{sequena}
_kV_n = \min(r > \, _kV_{n-1}, \  J_r \in \, _k\XX), \ n = 1, 2, \ldots, \  _kV_0 = 0.
\end{equation}

Now, let us define the random sequence,
\begin{equation}\label{sequenabo}
(_kJ_n, \, _kX_n) = \left\{
\begin{array}{ll}
(J_0, 0) & \ \text{for}  \ n = 0, \vspace{2mm} \\
(J_{\, _kV_n} \,,
\sum_{r = \, _kV_{n-1} +1}^{\, _kV_n} X_r) & \ \text{for}  \ n = 1, 2, \ldots.
\end{array}
\right.
\end{equation}

This sequence is also a Markov renewal process with phase space   $\XX \times [0, \infty)$, the initial distribution $\bar{p} = \langle p_i = \PP \{J_0 = i, X_0 = 0 \} = \PP \{J_0 = i \}, i \in {\mathbb X} \rangle$ and transition probabilities,
\begin{align}\label{numjabo}
_kQ_{ij}(t) & =  \PP \{ \, _kJ_{1} = j, \, _kX_{1} \leq t / \, _kJ_{0} = i, \, _kX_{0} = s \} \nonumber \\
& = Q_{ij}(t) + \sum_{n = 0}^\infty Q_{ik}(t) * Q^{(*n)}_{kk}(t) * Q_{kj}(t),  \,  t \geq 0, \, i, j \in \XX.
\end{align}

Here, symbol $*$ is used to denote the convolution of distribution functions (possibly improper),  and
$Q^{(*n)}_{kk}(t)$ is the $n$ times convolution of the distribution function  $Q_{kk}(t)$.

In this case, the Markov chain $_kJ_n$ has the transition probabilities,
\begin{align}\label{numjabonop}
_kp_{ij} =  \, _kQ_{ij}(\infty) & =  \PP \{ \, _kJ_{1} = j,  / \, _kJ_{0} = i \} \nonumber \\
& = p_{ij} + \sum_{n = 0}^\infty p_{ik} p_{kk}^n p_{kj} = p_{ij} + p_{ik} \frac{p_{kj} }{1 - p_{kk}}, \,  i, j \in \XX.
\end{align}

Note that condition ${\bf A}$ implies that  probabilities $p_{kk} \in [0, 1), \, k \in \XX$. 

The transition distributions for the Markov chain $_kJ_n$  are concentrated on the reduced phase space $_k\XX$, i.e., for every $i \in \XX$,
\begin{align}\label{sequenani}
\sum_{j \in \, _k\XX} \, _kp_{ij} & =  \sum_{j \in \, _k\XX} p_{ij} + p_{ik}  \sum_{j \in \, _k\XX} \frac{p_{kj} }{1 - p_{kk}} \nonumber \\
& =  \sum_{j \in \, _k\XX} p_{ij} + p_{ik}  = 1.
\end{align}

If the initial distribution $\bar{p}$ is concentrated on the phase space $_k\XX$, i.e., $p_k = 0$, then the random sequence  $(_kJ_{n}, \, _kX_{n}), n = 0, 1, \ldots$ 
can be considered as a Markov renewal process with the reduced  phase  $_k\XX \times [0, \infty)$, the initial distribution $_k\bar{p} = \langle \, p_i = \PP \{  _kJ_{0} = i, \, _kX_{0} = 0 \}  =  \PP \{  _kJ_{0} = i \} , i \in \, _k\XX \rangle$ and transition probabilities $_kQ_{ij}(t), t \geq 0, i, j \in \, _k\XX$.

If the initial distribution $\bar{p}$ is not concentrated on the phase space $_k\XX$, i.e., $p_k > 0$, then the random sequence  $(_kJ_{n}, \, _kX_{n}), n = 0, 1, \ldots$ 
can be interpreted as a Markov renewal process with so-called transition period.

Let us now introduce the  semi-Markov process, 
\begin{equation}\label{sepra}
_kJ(t) = \, _kJ_{_kN(t)}, \ t \geq 0,
\end{equation}
where $_kN(t) = \max(n \geq 0: \, _kT_n \leq t)$
is a number of jumps at time interval  $[0, t]$, for $t \geq 0$, and $_kT_n = \, _kX_1 + \cdots + \, _kX_n, \ n = 0, 1, \ldots$
are sequential moments of jumps,  for the semi-Markov process $_kJ(t)$.

As follows from the above remarks, the  semi-Markov process $_kJ(t), t \geq 0$ has transition probabilities $_kQ_{ij}(t), t \geq 0, i, j \in \, \XX$ 
concentrated on the reduced phase space  $_k\XX$, which can be interpreted as the actual  ``reduced'' phase space of this semi-Markov process $_kJ(t)$. 

If the initial distribution $\bar{p}$ is concentrated on the phase space $_k\XX$, then process $_kJ(t), t \geq 0$ can be considered as the semi-Markov process  
with the reduced  phase  $_k\XX$, the initial distribution $_k\bar{p} = \langle \, _kp_i = \PP \{_kJ_1(0) = i \}, i \in \, _k\XX \rangle$ and transition probabilities 
$_kQ_{ij}(t), t \geq 0,  i, j \in \, _k\XX$.

According to the above remarks, we can refer to the process $_kJ(t)$ as a reduced semi-Markov process. 

If the initial distribution $\bar{p}$ is not concentrated on the phase space $_k\XX$, then the process $_kJ(t), t \geq 0$ can be interpreted 
as a reduced semi-Markov process with transition period. \vspace{1mm}

{\bf 3.2. Transition characteristics  for reduced semi-Markov processes}. 
Relation (\ref{numjabonop}) implies the following formulas, for probabilities $_kp_{kj}$ and  $_kp_{ij}$, $i, j \in \, _k\XX$, 
\begin{equation}\label{inadopionnmoo}
\left\{
\begin{array}{lllll}
_kp_{kj} & =  \frac{p_{kj}}{1 - p_{kk}}, \\
_kp_{ij} & = p_{ij} + p_{ik} \, _kp_{kj} = p_{ij} + p_{ik} \frac{p_{kj}}{1 - p_{kk}}.
\end{array}
\right.
\end{equation}

It is useful to note that the second formula in relation (\ref{inadopionnmoo}) reduces to the first one, if to assign $i = k$ in this formula

Taking into account that $_kV_1$ is  Markov time for the  Markov renewal process $(J_n, X_n)$, we can write down  the  following system of stochastic equalities,  
for every $i, j \in \, _k\XX$,
\begin{equation}\label{inadopionxz}
\left\{
\begin{array}{lllll}
_kX_{i, 1} I( \, _kJ_{i, 1} = j) & \overset{d}{=} X_{i, 1} I(J_{i, 1} = j) \vspace{1mm} \\
& \ \ \ + \, (X_{i, 1} +  \, _kX_{k, 1}) I( J_{i, 1} = k) I(\, _kJ_{k, 1} = j),  \vspace{1mm} \\
_kX_{k, 1} I( \, _kJ_{k, 1} = j) &  \overset{d}{=} X_{k, 1} I(J_{k, 1} = j)   \vspace{1mm} \\
& \ \ \ + \, (X_{k, 1} +\, _kX_{k, 1}) I( J_{k, 1} = k) I(\, _kJ_{k, 1} = j),
\end{array}
\right.
\end{equation}
where: (a) $(J_{i, 1}, X_{i, 1})$ is a random vector, which  takes values in space ${\mathbb X} \times [0, \infty)$ and has the distribution $P \{J_{i, 1} = j, X_{i, 1} \leq t \} = Q_{ij}(t), \,  j \in \XX,  \, t \geq 0$, for every $i \in \XX$;
(b) $(_kJ_{i, 1}, \, _kX_{i, 1})$ is a random vector which takes values in the space $_k\mathbb{X} \times [0, \infty)$ and has distribution   $\mathsf{P} \{_kJ_{i, 1} = j, \, _kX_{i, 1} \leq t \} = \mathsf{P}_i \{ _kJ_{1} = j, \, _kX_{1} \leq t \}$
$= \, _kQ_{ij}(t), \, j \in \, _k\mathbb{X},  \, t \geq 0$, for every $i \in \XX$; (c) $(J_{i, 1}, X_{i, 1})$ and $(_kJ_{k, 1}, \, _kX_{k, 1})$ are independent random vectors, for every $i \in \XX$.

Let us denote,
\begin{equation}\label{inad}
_ke_{ij}^{(r)} =
\mathsf{E}_i \, _kX_{1}^r I( \, _kJ_{1} = j) = \int_0^\infty t^r \, _kQ_{ij}(dt), \  r = 0, 1,  \ldots, \ i, j \in \, _k\XX.
\end{equation}

Note that, 
\begin{equation}\label{inadasb}
_ke_{ij}^{(0)} = \, _kp_{ij},  \ i \in \XX,  j \in \, _k\XX.
\end{equation}

By taking  moments of the  order $r$  in stochastic relations (\ref{inadopionxz}) we get, for every $i, j \in \, _k\XX$,  the following system of linear equations  for the moments $_ke^{(r)}_{kj}$, $_ke^{(r)}_{ij}$ for $r = 1, \ldots, d$,
\begin{equation}\label{inadopionse}
\left\{
\begin{array}{lllll}
_ke^{(r)}_{kj} = e^{(r)}_{kj} + \sum_{l=0}^{r-1} \binom{r}{l}e^{(r-l)}_{kk} \, _ke^{(l)}_{kj} + p_{kk} \, _ke^{(r)}_{kj}, \vspace{3mm} \\
_ke^{(r)}_{ij} = e^{(r)}_{ij} + \sum_{l=0}^{r-1} \binom{r}{l}e^{(r-l)}_{ik} \, _ke^{(l)}_{kj} + p_{ik} \, _ke^{(r)}_{kj},
\end{array}
\right.
\end{equation}

Relation (\ref{inadopionse}) implies the following recurrent formulas for moments  $_ke^{(r)}_{kj}$ and  $_ke^{(r)}_{ij}$, which should be used, for every $i, j \in \, _k\XX$, 
recurrently for $r = 1, \ldots, d$,  
\begin{equation}\label{inadopionmi}
\left\{
\begin{array}{lllll}
_ke^{(r)}_{kj}  & = \frac{1}{1 - p_{kk}} \big( e^{(r)}_{kj}
+ \sum_{l=0}^{r-1} \binom{r}{l} e^{(r-l)}_{kk} \, _ke^{(l)}_{kj} \big), \vspace{3mm} \\
_ke^{(r)}_{ij} & = e^{(r)}_{ij}+ \sum_{l=0}^{r-1} \binom{r}{l} e^{(r-l)}_{ik} \, _ke^{(l)}_{kj}   \vspace{3mm}  \\
& \quad +  \frac{p_{ik}}{1 - p_{kk}} \big( e^{(r)}_{kj}+ \sum_{l=0}^{r-1} \binom{r}{l} e^{(r-l)}_{kk} \, _ke^{(l)}_{kj} \big),
\end{array}
\right.
\end{equation}

It is useful to note that the second formula in relation (\ref{inadopionmi}) reduces to the first one, if to assign $i = k$ in this formula. \vspace{1mm}

{\bf 3.3. Hitting times for  reduced semi-Markov processes}. Let us assume that $k \neq 0$ and  introduce the 
first hitting time to state $0$ for the reduced semi-Markov process  $_kJ(t)$,
\begin{equation}\label{hittake}
_kW_0 = \inf(t \geq  \, _kX_1: \, _kJ(t) = 0) = \sum_{n = 1}^{_kU_0} \, _kX_n,
\end{equation}
where $_kU_0 =  \min(n \geq 1: \, _kJ_n = 0)$ is the first hitting time to state $0$ by the reduced  Markov chain $_kJ_n$.

Let also introduce moments,
\begin{equation}\label{mhittaka}
_kE_{i0}^{(r)} = \EE_i \, _kW^r_0,   \  r  = 0, 1, \ldots, d, \,  i \in \XX.
\end{equation}

Note that, 
\begin{equation}\label{mhittakani}
_kE_{i0}^{(0)} = 1, i \in {\mathbb X}.
\end{equation}

The following theorem plays the key role in what follows.
\vspace{1mm}

{\bf Theorem 1}. {\em Conditions ${\bf A}$, ${\bf B}$ and  ${\bf C_d}$ assumed to hold for the semi-Markov process $J(t)$ also hold for the reduced 
semi-Markov process  $_kJ(t)$, for any state $k \neq 0$. Moreover,   the hitting times $W_0$ and $_kW_0$ to the state $0$, respectively, for semi-Markov processes
$J(t)$ and $_kJ(t)$, coincide, and, thus, for every $r = 1, \ldots, d$ and $i \in \XX$,}
\begin{equation}\label{coin}
E_{i0}^{(r)} = \EE_i W_0^r = \, _kE_{i0}^{(r)}
= \EE_i \, _kW_0^r.
\end{equation}

{\bf Proof}. Holding of conditions ${\bf A}$ and ${\bf B}$ for  the semi-Markov process $_kJ(t)$ is obvious. Holding of condition  ${\bf C_d}$ for  the semi-Markov process $_kJ(t)$  
 follows from relation (\ref{inadopionmi}).

The first hitting times to a state $0$ are connected for Markov chains  $J_n$  and
 $_kJ_n$ by the following relation,
 \begin{equation}\label{relanaka}
 U_0  = \min(n \geq 1: J_n = 0)  = \min(_kV_n \geq 1: \, _kJ_n  = j) = \, _kV_{_kU_0},
 \end{equation}
where $_kU_0 = \min(n \geq 1: \, _kJ_n = 0)$.

The above relations imply that the following relation holds for the first hitting times to state $0$, for the semi-Markov
processes $J(t)$ and $_kJ(t)$,
\begin{equation}\label{relanak}
W_0 = \sum_{n = 1}^{U_0} X_n  = \sum_{n = 1}^{_kV_{_kU_0 }} X_n
=  \sum_{n = 1}^{_kU_0 } \, _kX_n = \, _kW_0.
\end{equation}

The equality for moments of the first hitting times  is an obvious corollary of relation (\ref{relanak}). $\Box$

We can  write down the recurrent systems of linear equations (\ref{inadopionba}) for moments $_kE_{k0}^{(r)}$  and  $_kE_{i0}^{(r)}, i \in \, _k\XX$ of the reduced semi-Markov 
process $_kJ(t)$, which should be solved recurrently, for $r = 1, \ldots, d$,
\begin{equation}\label{inadopionball}
\left\{
\begin{array}{ll}
_kE^{(r)}_{k0} & = \, _kf^{(r)}_{k0}   + \sum_{j \in \, _k\XX, j \neq 0} \, _kp_{kj} \, _kE^{(r)}_{j0}, \vspace{2mm} \\
_kE^{(r)}_{i0} & = \, _kf^{(r)}_{i0}   + \sum_{j \in \, _k\XX, j \neq 0} \, _kp_{ij} \, _kE^{(r)}_{j0},  i \in \, _k\XX,
\end{array}
\right.
\end{equation}
where
\begin{equation}\label{inadopnn}
_kf^{(r)}_{i0}  = \, _ke^{(r)}_{i0} + \sum_{j \in \, _k\XX, \, j \neq 0} \sum_{l=0}^{r-1} \binom{r}{l} \, _ke^{(r-l)}_{ij} \, _kE^{(l)}_{j0}, \ i \in \XX.
 \end{equation}

Theorem 1 makes it possible to compute moments $E_{i0}^{(r)} = \, _kE_{i0}^{(r)}, i \in \XX, r = 1, \ldots, d$ in the way alternative to
solving recurrent systems of linear equations (\ref{inadopionba}).

Instead of this, we can, first, compute transition probabilities and moments of transition  times for the reduced semi-Markov process $_kJ(t)$
using, respectively, relations (\ref{inadopionnmoo}) and  (\ref{inadopionmi}), and, then, by solving the systems of linear equations
(\ref{inadopionball}) sequentially for $r = 1, \ldots, d$.

Note that every system given in  (\ref{inadopionba}) has $m$ equations for moments  
$E^{(r)}_{i0}, i \in \, \XX, i \neq 0$ plus the explicit formula for computing moment $E^{(r)}_{00}$ as function of moments $E^{(r)}_{i0}, i \in \, \XX, i \neq 0$. 

While, every  system given in  (\ref{inadopionball}) has, in fact,  $m - 1$ equations for moments  
$_kE^{(r)}_{i0}, i \in \, _k\XX, i \neq 0$, plus two explicit formulas for computing moment $_kE^{(r)}_{00}$ and $_kE^{(r)}_{k0}$ as functions of moments $_kE^{(r)}_{i0}, i \in \, _k\XX, i \neq 0$. \\

{\bf 4, Algorithms of sequential phase space reduction} \\
 
In this section, we present a multi-step algorithm for sequential reduction of phase space for semi-Markov processes. We also present the recurrent algorithm for computing power moments of hitting times for semi-Markov processes, which are based on the above algorithm of sequential reduction 
of the phase space. \vspace{1mm}

{\bf 4.1. Sequential reduction of phases space for semi-Markov processes}. In what follows,  let $i \in \{1, \ldots, m \}$ and let  $\bar{k}_{i, m}  = \langle k_{i, 1}, \ldots, k_{i, m} \rangle =  \langle k_{i, 1}, \ldots$, $k_{i, m -1}, i \rangle$ be a permutation of the sequence  $\langle 1, \ldots, m \rangle$ such that $k_{i, m} = i$, and let $\bar{k}_{i, n} = \langle  k_{i, 1}, \ldots, k_{i, n} \rangle$, $n = 1, \ldots, m$ be the corresponding chain of growing sequences of states from space $\XX$.

 Let us assume that $p_0 + p_i = 1$. Denote as $_{\bar{k}_{i, 0}}J(t) = J(t)$, the initial semi-Markov process. Let us exclude  state  $k_{i, 1}$ from the phase space $_{\bar{k}_{i, 0}}\XX = \XX$ of semi-Markov process $_{\bar{k}_{i, 0}}J(t)$ using the time-space screening procedure described in Section 3. Let
$_{\bar{k}_{i, 1}}J(t)$ be the corresponding reduced semi-Markov process. The above procedure can be repeated. The state $k_{i, 2}$ can  be excluded from the phase space of the semi-Markov process $_{\bar{k}_{i, 1}}J(t)$. Let $_{\bar{k}_{i, 2}}J(t)$ be the corresponding reduced semi-Markov process. By continuing the above procedure for states $k_{i, 3}, \ldots, k_{i, n}$,  we construct the reduced semi-Markov process $_{\bar{k}_{i, n}}J(t)$.

The process  $_{\bar{k}_{i, n}}J(t)$ has,  for every $n = 1, \ldots, m$,   the actual  ``reduced''  phase space, 
\begin{equation}\label{inoonan}
_{\bar{k}_{i, n}}\XX = \, _{\bar{k}_{i, n-1}}\XX \setminus \{ k_{i, n} \} =   \XX \setminus \{ k_{i, 1}, k_{i, 2}, \ldots, k_{i, n} \}.
\end{equation}

The transition probabilities  $_{\bar{k}_{i, n}}p_{k_{i, n}, j'}$, $_{\bar{k}_{i, n}}p_{i'j'}$, $i', j' \in \, _{\bar{k}_n}\XX$, and the moments $_{\bar{k}_{i, n}}e_{k_{i, n}, j'}^{(r)}$,   
$_{\bar{k}_{i, n}}e_{i'j'}^{(r)}, i', j' \in \, 
_{\bar{k}_{i, n}}\XX, \, r = 1, \ldots, d$ are determined for the semi-Markov process $_{\bar{k}_{i, n}}J(t)$ by the transition probabilities and the expectations of sojourn  times for the semi-Markov process $_{\bar{k}_{i, n-1}}J(t)$, respectively, via relations (\ref{inadopionnmoo}) and (\ref{inadopionmi}), which take the following  recurrent forms, for $i', j' \in \, 
_{\bar{k}_{i, n}}\XX, r = 1, \ldots, d$ and $n = 1, \ldots, m$,  
\begin{equation}\label{innmoonan}
\left\{
\begin{array}{lllll}
_{\bar{k}_{i, n}}p_{k_{i, n}, j'} & = \frac{_{\bar{k}_{i, n-1}}p_{k_{i, n}, j'}}{1 - \, _{\bar{k}_{i, n-1}}p_{k_{i, n}, k_{i, n}}}, \vspace{2mm} \\
_{\bar{k}_{i, n}}p_{i' j'} & = \,  _{\bar{k}_{i, n-1}}p_{i' j'} \vspace{2mm}  \\ 
& \quad + \, _{\bar{k}_{i, n-1}}p_{i' k_{i, n}} \frac{_{\bar{k}_{i, n-1}}p_{k_{i, n}, j'}}{1 - \, _{\bar{k}_{i, n-1}}p_{k_{i, n}, k_{i, n},}}, 
\end{array}
\right.
\end{equation}
and
\begin{equation}\label{dopionmina}
\left\{
\begin{array}{lllll}
_{\bar{k}_{i, n}}e^{(r)}_{ k_{i, n} j'}  & = \frac{1}{1 - \, _{\bar{k}_{i, n-1}}p_{k_{i, n} k_{i, n}}} \big( e^{(r)}_{k_{i, n} j'} \vspace{2mm}  \\ 
& \quad + \sum_{l=0}^{r-l} \binom{r}{l} \, _{\bar{k}_{i, n-1}}e^{(r-l)}_{k_{i, n} k_{i, n}} \, _{\bar{k}_{i, n-1}}e^{(l)}_{k_{i, n} j'} \big), \vspace{2mm} \\
\, _{\bar{k}_{i, n}}e^{(r)}_{i' j'} & = \, _{\bar{k}_{i, n-1}}e^{(r)}_{i' j'}+ \sum_{l=0}^{r-1} \binom{r}{l} \, _{\bar{k}_{i, n-1}}e^{(r-l)}_{i' k_{i, n} } \,  _{\bar{k}_{i, n-1}}e^{(l)}_{k_{i, n}  j'}   \vspace{2mm}  \\
& \quad +  \frac{_{\bar{k}_{i, n-1}}p_{i' k_{i, n} }}{1 - \, _{\bar{k}_{i, n-1}}p_{k_{i, n} k_{i, n} }} \big( \, _{\bar{k}_{i, n-1}}e^{(r)}_{k_{i, n}  j'} \vspace{2mm}  \\
& \quad + \sum_{l=0}^{r-l} \binom{r}{l} \, _{\bar{k}_{i, n-1}}e^{(r-l)}_{k_{i, n} k_{i, n} } \, \, _{\bar{k}_{i, n-1}}e^{(l)}_{k_{i, n} j'} \big).
\end{array}
\right.
\end{equation}
\vspace{1mm}

{\bf 4.2. Recurrent algorithms for computing of moments of hitting times}. Let us $_{\bar{k}_{i, n}}W_{0}$ be  the first hitting time to state $0$ for the reduced semi-Markov process  $_{\bar{k}_{i, n}}J(t)$ and $_{\bar{k}_{i, n}}E_{i' 0}^{(r)} = \EE_{i'} \, _{\bar{k}_{i, n}}W_{0}^r,  i' \in \, _{\bar{k}_{i, n}}\XX, r = 1, \ldots, d$ be the  moments for these random variables.

By Theorem 1, the above moments of  hitting time coincide for the semi-Markov processes  $_{\bar{k}_{i, 0}}J(t)$,  $_{\bar{k}_{i, 1}}J(t), \ldots, \, _{\bar{k}_{i, n}}J(t)$, i.e., for $n' = 0, \ldots, n$,
\begin{equation}\label{vopty}
_{\bar{k}_{j, n'}}E_{k_{i, n'} 0}^{(r)} = E_{k_{i, n'} 0}^{(r)}, \ _{\bar{k}_{j, n'}}E_{i' 0}^{(r)} = E_{i' 0}^{(r)}, \, i' \in \, _{\bar{k}_{i, n}}\XX, \, r = 1, \ldots, d.
\end{equation}

Moreover, the moments of  hitting times $_{\bar{k}_{j, n}}E_{k_{i, n} 0}^{(r)}$, $_{\bar{k}_{i, n}}E_{i' 0}^{(r)}, i' \in \, _{\bar{k}_{i, n}}\XX, r = 1, \ldots, d$   resulted by  the recurrent algorithm of sequential phase space reduction described above, are invariant with respect to any permutation $\bar{k}'_{i, n} = \langle k'_{i, 1}, \ldots$, $k'_{i, n} \rangle$ of sequence $\bar{k}_{i, n} =  \langle k_{i, 1}, \ldots, k_{i, n}  \rangle$.

Indeed, for every permutation $\bar{k}'_{i, n}$ of sequence $\bar{k}_{i, n}$, the corresponding reduced semi-Markov process $_{\bar{k}'_{i, n}}J(t)$ is constructed  as the sequence of states for the initial semi-Markov process $J(t)$   at sequential moment of its hitting into the same reduced phase space $_{\bar{k}'_{i, n}}\XX = \XX \setminus \{ k'_{i, 1}, \ldots, k'_{i, n} \} = \, _{\bar{k}_{i, n}}\XX  = \XX \setminus \{ k_{i, 1}, \ldots, k_{i, n} \}$. The times between sequential jumps  of the  reduced semi-Markov process $_{\bar{k}'_{i, n}}J(t)$ are the  times between sequential hitting of the above  reduced phase space by the initial  semi-Markov process $J(t)$.

This implies that the transition probabilities  $_{\bar{k}_{i, n}}p_{k_{i, n} j'}$, $_{\bar{k}_{i, n}}p_{i' j'}, i',  j' \in \, _{\bar{k}_{i, n}}\XX$ and the moments  $_{\bar{k}_{i, n}}e_{k_{i, n} j'}^{(r)}$, 
$_{\bar{k}_{i, n}}e_{i'j'}^{(r)}, i', j' \in \, _{\bar{k}_{i, n}}\XX, r = 1, \ldots, d$ and, in sequel,  moments $_{\bar{k}_{i, n}}E_{k_{i, n} 0}^{(r)}$, $_{\bar{k}_{i, n}}E_{i' 0}^{(r)}, i' \in \, _{\bar{k}_{i, n}}\XX, r = 1, \ldots, d$  are, for every $n = 1, \ldots, m$,  invariant  with respect to any permutation $\bar{k}'_{i, n}$ of the sequence $\bar{k}_{i, n}$.

Let us now choose $n = m$. In this case, the reduced semi-Markov process $_{\bar{k}_{i, m}}J(t)$ has the one-state phase space $_{\bar{k}_{i, m}}\XX = \{ 0 \}$ and state $k_{i, m} = i$. 

In this case, the reduced semi-Markov process $_{\bar{k}_{i, m}}J(t)$ return to state $0$ after every jump and hitting time to state $0$ coincides with the sojourn time in state $_{\bar{k}_{i, m}}J(0)$.

Thus, the transition probabilities,
\begin{equation}\label{inadopiobaha}
_{\bar{k}_{i, m}}p_{i 0}  = \, _{\bar{k}_{i, m}}p_{0 0} = 1.
 \end{equation} 
 
 Also, by Theorem 1, moments,
\begin{equation}\label{inadopiobahar}
E^{(r)}_{i 0} = \, _{\bar{k}_{i, m}}E^{(r)}_{i 0} = \, _{\bar{k}_{i, m}}e^{(r)}_{ i 0}, r = 1, \ldots, d,
\end{equation}
 and 
\begin{equation}\label{inadopiobahara}
E^{(r)}_{0 0}  =  \, _{\bar{k}_{i, m}}E^{(r)}_{0 0} 
= \, _{\bar{k}_{i, m}}e^{(r)}_{ 0 0}, r = 1, \ldots, d.
\end{equation}

The above remarks can be summarized in the following theorem, which presents the recurrent algorithm for computing of power moments for  hitting times. 
\vspace{1mm}

{\bf Theorem 2}. {\em Moments $E^{(r)}_{i 0}, E^{(r)}_{0 0}, r = 1, \ldots, d$ are given, for every $i = 1, \ldots, m$, by 
formulas  {\rm (\ref{inadopiobahar}) -- (\ref{inadopiobahara})}, where transition probabilities $_{\bar{k}_{i, n}}p_{k_{i, n}, j'}$, $_{\bar{k}_{i, n}}p_{i'j'}$, $i', j' \in \, _{\bar{k}_n}\XX$, and moments $_{\bar{k}_{i, n}}e_{k_{i, n}, j'}^{(r)}$,   
$_{\bar{k}_{i, n}}e_{i'j'}^{(r)}, i', j' \in \,  _{\bar{k}_{i, n}}\XX, \, r = 1, \ldots, d$  are determined, for $n = 1, \ldots, m$,
by recurrent formulas {\rm (\ref{innmoonan}) -- (\ref{dopionmina})} and formula  {\rm (\ref{inadopiobaha})}. The moments 
$E^{(r)}_{i 0}, E^{(r)}_{0 0}, r = 1, \ldots, d$  are invariant with respect to any permutation $\bar{k}_{i, m}$ of 
sequence   $\langle 1, \ldots, m \rangle$ used in the above recurrent algorithm.}  \\

{\bf 5. Generalizations and examples} \\

In this section, we describe several  variants for generalization of the results concerned recurrent algorithms for computing power moments of hitting times and accumulated rewards of hitting type. \vspace{1mm}

{\bf 5.1. Real-valued accumulated rewards of hitting type}. First, we would like to mention that Theorems 1 and 2 can be  generalized on the model, where of the Markov renewal process $(J_n, X_n), n = 0, 1, \ldots$ has the phase space 
${\mathbb X} \times \mathbb{R}_1$, an initial distribution $\bar{p} = \langle p_i = \PP \{J_0 = i, X_0 = 0 \} = \PP \{J_0 = i \}, i \in {\mathbb X} \rangle$ and transition probabilities,
\begin{equation}\label{semikav}
Q_{ij}(t) = \PP \{ J_{1} = j, X_{1} \leq t / J_{0} = i, X_{0} = s \}, \ (i, s), (j, t) \in  {\mathbb X} \times  \mathbb{R}_1.
\end{equation}

In this case, we the random variable,
\begin{equation}\label{semikavny}
W_0 = \sum_{n = 1}^{U_0} X_n
\end{equation} 
can  be  be interpreted as a reward accumulated  on trajectories of Markov chain $J_n$ 
up to its first hitting time $U_0 = \min(n \geq 1, J_n = 0)$ of this Markov chain  to the state $0$.

Condition  ${\bf C_d}$ should be replaced by condition:
\begin{itemize}
\item[${\bf \dot{C}_d}$:]  $\EE_i |X_1|^d  < \infty, i \in \mathbb{X}$.
\end{itemize}

As well known, in this case moments $\dot{E}_{i}^{(d)} = \EE_i |W_0|^d, i \in \XX$ are finite.

All recurrent relations for moments  $\EE_{i}^{(r)} = \EE_i W_0^r, r = 1, \ldots, d, i \in \XX$, given in Sections 3 -- 4,  as well as Theorems 1 and 2 take the same forms as in the case of nonnegative rewards. \vspace{1mm}

{\bf 5.2. Vector accumulated rewards of hitting type}. Second, we would like to show, how the above results can be generalized on the case of vector accumulated rewards. 

For simplicity, let us consider the bivariate case, where the Markov renewal process $(J_n, \bar{X}_n) = (J_n , (X_{1, n}, X_{2, n})) = 0, 1, \ldots$ has the phase space 
${\mathbb X} \times \mathbb{R}_2$, an initial distribution $\bar{p} = \langle p_i = \PP \{J_0 = i, \vec{X}_0 = (0, 0) \} =  \PP \{J_0 = i \}, i \in {\mathbb X} \rangle$ and transition probabilities,
\begin{equation}\label{semikava}
Q_{ij}(\bar{t}) = \PP \{ J_{1} = j, \bar{X}_{1} \leq \bar{t} / J_{0} = i, \bar{X}_{0} = s \}, \ (i, \bar{s}), (j, \bar{t}) \in  {\mathbb X} \times  \mathbb{R}_2. 
\end{equation}

Here and henceforth symbol $\bar{u} \leq \bar{v}$ for vectors $\bar{u} = (u_1, u_2), \bar{v} = (v_1, v_2) \in \mathbb{R}_2$ means that $u_1 \leq v_1, u_2 \leq v_2$.

The vector accumulated reward $\bar{W}_0 = (W_{1, 0}, W_{2, 0})$ is defined as a bivariate random vector with components,
\begin{equation}\label{seava}
W_{l, 0} = \sum_{n = 1}^{U_0} X_{l, n}, \, l = 1, 2.
\end{equation}

Condition ${\bf \dot{C}_d}$ should be replaced by condition:
\begin{itemize}
\item[${\bf \dot{C}'_d}$:]  $\EE_i |X_{l, 1}|^d < \infty, l= 1, 2, \, i \in \mathbb{X}$.
\end{itemize}

In this case, moments $\dot{E}_{l, i}^{(d)} = \EE_i |W_{l, 0}|^d < \infty,\,  l = 1, 2, \, i \in \XX$. 

Let us introduce mixed moments,
\begin{equation}\label{seavamu}
E_i^{(q, r)} = \EE_{i} W^{q}_{1, 0}W^{r - q}_{2, 0}, \, 0 \leq q \leq r \leq d, \, i \in \XX.
\end{equation}

Let us define random variables $W_0(a) = a W_{1, 0}  + (1-a)W_{2,  0}, 0 \leq a \leq 1$. By definition, $W_0(a) =
\sum_{n =1}^{U_0} (a X_{1, n} + (1-a) X_{2, n})$ is also an accumulated reward for the corresponding  local rewards 
$X_n(a) = a X_{1, n} + (1-a) X_{2, n}, n = 1, 2, \ldots$.

Let us denote $E_i^{(r)}(a) = \EE_{i} W^r(a), 0 \leq a \leq 1, r = 1, \ldots, d, i \in \XX$. 

Let us assume that the moments of non-negative accumulated rewards   $E_i^{(r)}(a)$ are found for $r+1$ values $a_p, p = 0, \ldots, r$, for example, for $a_p = \frac{p}{r}, p = 0, \ldots, r$, for every  $r = 1, \ldots, d, i \in \XX$ using recurrent algorithms described in Sections 2 -- 4.

Then,  the following system of linear equations can be written down for the correlation moments
$E_i^{(q, r)}, q = 0, \ldots, r$, for every $r = 1, \ldots, d, i \in \XX$, 
\begin{equation}\label{ikavany}
\Big\{ E_i^{(r)}(a_p) = \sum_{q = 0}^r \binom{r}{q} a_p^{q}(1-a_p)^{r - q} E_i^{(q, r)},  \ q = 0, \ldots, r.
\end{equation}

It can be shown that the above linear system  has the non-zero determinant. Thus, the moments $E_i^{(r)}(a_p)$, $a_p = \frac{p}{r}$, $p = 0, \ldots, r$ uniquely determine the mixed  moments $E_i^{(q, r)}, 0 \leq q \leq r$, for every $r = 1, \ldots,  d, \,  i \in \XX$. \vspace{1mm}

{\bf 5.3. General hitting times with hitting state indicators}. Third, the above results can be generalized on the case of more general hitting times, 
\begin{equation}\label{mhivny}
W_{\DD} = \sum_{n = 1}^{U_{\DD}} X_n,
\end{equation} 
where $U_{\DD} = \min(n \geq 1, J_n \in  \DD)$, for some nonempty set $\DD \subset \XX$.

In this case main object of studies are  power moments for the first hitting times with hitting state indicators,
\begin{equation}\label{mhittany}
E_{\DD, ij}^{(r)} = \EE_i W_{\DD}^r I( J_{U_{\DD}} = j),   \  r  = 0, 1, \ldots, d, \,   j \in \DD,  i \in \XX.
\end{equation}

Note that,
\begin{equation}\label{mhittanyny}
E_{\DD, ij}^{(0)}  = \PP_i \{J_{U_{\DD}} = j \}, i \in {\mathbb X}, j \in \DD.
\end{equation}

As well known, conditions ${\bf A}$,  ${\bf B}$ and ${\bf C_d}$ imply that, for any nonempty set  $\DD \subset \XX$,
\begin{equation}\label{mhittavny}
E_{\DD, ij}^{(r)}  < \infty, \  r = 1, \ldots, d,  \ i \in \XX, j \in \DD.
\end{equation}

Note that the simpler condition  ${\bf A}$ can, in fact, be replaced by a simpler condition: 
\begin{itemize}
\item[${\bf  A_{\DD}}$:] $\PP_i \{ U_{\DD}   < \infty \} = 1, \, i \in \XX$.
\end{itemize}

In this case,  theorems, analogous of Theorems 1 and 2, take place,  and  recurrent systems of linear equations and recurrent formulas analogous to those given in Sections 2 -- 4
can  be written down.   

For example, let $_kE_{\DD, ij}^{(r)},  r = 1, \ldots, d,  \ i \in \XX, j \in \DD$ be the moments $E_{\DD, ij}^{(r)}  < \infty, \  r = 1, \ldots, d,  \ i \in \XX, j \in \DD$ computed for the reduced semi-Markov process $_kJ(t)$, for some $k  \notin \DD$.

The key recurrent systems of linear equations analogous to (\ref{inadopionball})
take, for every $ j \in \DD$, nonempty set  $\DD \subset \XX$ and $k \notin \DD$, the following form, for $r = 0, \ldots, d$,
\begin{equation}\label{inadopionballbb}
\left\{
\begin{array}{ll}
_kE_{\DD, kj}^{(r)}  & = \, _kf^{(r)}_{\DD, kj}   + \sum_{j' \in \, _k\XX \setminus \DD} \, _kp_{kj'} \, kE_{\DD, j' j}^{(r)}, 
\vspace{2mm} \\
_kE_{\DD, ij}^{(r)}  & = \, _kf^{(r)}_{\DD, ij}    + \sum_{j' \in \, _k\XX \setminus \DD} \, _kp_{ij'} \, _kE_{\DD, j' j}^{(r)},  
\, i \in \, _k\XX,
\end{array}
\right.
\end{equation}
where
\begin{equation}\label{inadopnnbb }
_kf^{(r)}_{\DD, ij}   = \, _ke^{(r)}_{ij} + \sum_{j' \in \, _k\XX \setminus \DD} 
\sum_{l=0}^{r-1} \binom{r}{l} \, _ke^{(r-l)}_{ij'} \, _kE_{\DD, j' j}^{(r)}, \ i \in \XX.
 \end{equation}

The difference with the recurrent systems of linear equations (\ref{inadopionball}) is that, in this case, the 
corresponding system of linear equations for hitting probabilities $E_{\DD, ij}^{(0)}, i \in \XX$ should also be  solved.

Also, the corresponding changes caused by  replacement of the hitting state $0$ by state $j \in \DD$ and set $_k\XX \setminus \{ 0 \}$ by set  $_k\XX \setminus \DD$ sould be taken into account when writing down  systems of  linear equations (\ref{inadopionballbb}) instead of systems of  
linear equations  (\ref{inadopionball}). \vspace{1mm}

{\bf 5.4. Place-dependent hitting times}. Fourth, the above results can be generalized on so-called place-dependent hitting times,
\begin{equation}\label{inanbb} 
Y_{\GG} = \sum_{n = 1}^{U_{\GG}} X_n, 
\end{equation}
where $U_{\GG} = \min(n \geq 1: (J_{n-1}, J_n) \in  \GG)$, for some nonempty 
set $\GG \subset \XX \times \XX$.

Note that set $\GG$  can be represented in the form $\GG = \cup_{i \in \XX} \,  \{ i \} \times \GG_i$, where $\GG_i = \{ j \in \XX: (i, j) \in \GG \}$. Respectively, 
the first hitting time $U_{\GG}$ can be represented as $U_{\GG} = \min(n \geq 1: J_n\in  \GG_{J_{n-1}})$. This representation explains using of the term ``place-dependent hitting time''. 

In fact, the above model can be embedded in the previous one, if to consider the new Markov renewal process $(\bar{J_n}, X_n) = 
((J_{n -1}, J_n), X_n), n = 0,1, \ldots$ constructed from the initial Markov renewal process $(J_n, X_n), n = 0, 1, \ldots$ by aggregating sequential states 
for the initial embedded Markov chain $J_n$.

The Markov renewal process  $(\bar{J_n}, X_n)$ has the phase space $(\XX \times \XX) \times [0, \infty)$. For simplicity, we can take the initial state 
$\bar{J_0} = (J_{-1}, J_0)$, where $J_{-1}$ is a random variable taking values in space $\XX$ and independent on the Markov renewal process $(J_n, X_n)$.

Note that the simpler condition  ${\bf A}$ can, in fact, be replaced by a simpler condition: 
\begin{itemize}
\item[${\bf  A'_{\GG}}$:] $\PP_i \{ U_{\GG}   < \infty \} = 1, \, i \in \XX$.
\end{itemize}

The above assumption, that domain $\GG$ is hittable,  is implied by condition  ${\bf A}$,  for any domain $\GG$ containing a pair of states $(i, j)$ such that  $p_{ij} > 0$. 

The results  concerned moments of usual accumulated rewards $W_{\DD}$  can be expanded to the place-depended  accumulated rewards $\YY_{\GG}$ for hittable domains,  using the above embedding procedure.  \vspace{1mm}

{\bf 5.5. Time-dependent hitting times}. Let   $(J_n, X_n), n = 0, 1, \ldots$ be an  inhomogeneous in time Markov renewal process, i.e., an  inhomogeneous in time  Markov chain with phase space with the phase space ${\mathbb X} \times [0, \infty)$, an initial distribution $\bar{p} = \langle p_i = \PP \{J_0 = i, X_0 = 0 \} = \PP \{J_0 = i \}, i \in {\mathbb X} \rangle$ and transition probabilities, defined for $(i, s), (j, t) \in  {\mathbb X} \times [0, \infty)$ and $n = 0, 1, 2, \ldots$,
\begin{equation}\label{semikad}
Q^{(n+1)}_{ij}(t) = \PP \{ J_{n+1} = j, X_{n+1} \leq t / J_{n} = i, X_{n} = s \}.
\end{equation}

As in homogeneous in time case, we exclude instant jumps and assume that the following condition holds;
\begin{itemize}
\item[${\bf B'}$:]  $Q^{(n)}_{ij}(0) = 0, \, i, j \in \mathbb{X}, \, n \geq 1$.
\end{itemize}

Process $(J_n, X_n)$ can be transformed in a homogeneous in time Markov renewal process by adding  to  this process an additional counting time component $J'_n = n, n = 0, 1, \ldots$. Indeed, process  $(\bar{J}_n, X_n) = ((J'_n, J_n), X_n), n = 0, 1, \ldots$ is a homogeneous in time  Markov renewal process. This process has the phase space $(\NN \times \XX) \times [0, \infty)$, where $\NN = \{ 0, 1, \ldots \}$. It has the  initial distribution $\bar{p} = \langle p_i = \PP \{J'_0 = 0, J_0 = i, X_0 = 0 \} =  \PP \{J_0 = i \}, i \in {\mathbb X} \rangle$
and  transition probabilities,
{\small
\begin{equation}\label{semikadbb}
Q_{(n,i),(k,j)}(t)  = \left\{
\begin{array}{lll}
Q^{(n+1)}_{ij}(t) & \text{for} \ t \geq 0, k = n + 1, n = 0, 1, \ldots, i, j \in \XX, \vspace{2mm} \\
\quad 0 & \text{for} \ t \geq 0, k \neq  n + 1, n = 0, 1, \ldots, i, j \in \XX.  
\end{array}
\right.
\end{equation}
}

The phase space of the process $(\bar{J}_n, X_n)$ is countable. 

Let now define a time-truncated version of process $(\bar{J}_n, X_n)$ as the process $(\bar{J}^{(h)}_n, X^{(h)}_n) = ((J'_{n \wedge h}, J_{n \wedge h}), X_{n \wedge h}), n = 0, 1, \ldots$, for some integer $h \geq 1$. 

The process  $(\bar{J}^{(h)}_n, X^{(h)}_n), n = 0, 1, \ldots$ is also a homogeneous in time Markov renewal process. It has the finite phase space $(\HH \times \XX) \times [0, \infty)$, where 
$\HH =  \{ 0, 1, \ldots, h \}$.

Let $ \langle \DD_1, \ldots, \DD_h \rangle$ be some sequence of subsets of space $\XX$ such that  $\DD_h = \XX$ and let $U_{\tilde{\DD}_h} =  \min( n \geq 1: \bar{J}^{(h)}_n \in \{n \} \times \DD_n) =  \min( n \geq 1: J_n \in \DD_n)$ is the first hitting time to the domain $\tilde{\DD}_h = \cup_{n = 1}^h  \{n \} \times \DD_n$ for the Markov
chain $ \bar{J}^{(h)}_n$. 

Obviously, $\PP_i \{ U_{\tilde{\DD}_h} \leq h \} = 1, i \in \XX$, i.e., domain $\tilde{\DD}$ is hittable for the Markov chain $\bar{J}^{(h)}_n$. 

Thus, all results presented in Sections 2 -- 4 can be applied to the time-dependent accumulated rewards of hitting type,
\begin{equation}\label{inanbbol} 
Z_{\tilde{\DD}_h} = \sum_{n = 1}^{U_{\tilde{\DD}_h}} X_n.  
\end{equation}

Note only hat condition ${\bf C_d}$ should be, in this case, replaced by condition:
\begin{itemize}
\item[${\bf C_{h, d}}$:]  $\EE \{ X_{n}^d I(J_n = j) / J_{n -1} = i \}  < \infty, n = 1, \ldots, h,  i, j \in \mathbb{X}$.
\end{itemize}

In conclusion,  we would like also to note that it is possible to combine all five listed above generalization aspects 
in the frame of one semi-Markov model. \vspace{1mm}

{\bf 5.6. An example}. Let us  consider a numerical example illustrating the recurrent algorithm for computing power moment of hitting times and accumulated rewards of hitting times for semi-Markov processes, based on sequential  reduction of their phase spaces.

Let $J(t)$  be a semi-Markov process 
with the phase space $\XX = \{0, 1, 2, 3 \}$, and  the $4 \times 4$ matrix  of transition probabilities, $\| Q_{ij}(t) \|$, which has the following form, for $t \geq 0$,
\begin{equation}\label{partysama}
\left \| 
\begin{array}{lllllll}
0 & 0  & 0 & I(t \geq 1) \vspace{1mm} \\
\frac{1}{4}(1 - e^{-t/4})  & \frac{1}{4}(1 - e^{-t/4})   & \frac{1}{4}(1 - e^{-t/4})     &  \frac{1}{4}(1 - e^{-t/4}) \vspace{1mm} \\
0 &  \frac{1}{3}(1 - e^{- t/3})  & \frac{1}{3}(1 -  e^{- t/3})  &  \frac{1}{3}(1 -  e^{- t/3})  \vspace{1mm} \\
0 &  0 & \frac{1}{2}I(t \geq 2)  &  \frac{1}{2}I(t \geq 2)  \\
\end{array}
\right \|.
\end{equation} 

The  $4 \times 4$ matrices  of transition probabilities $\| p_{ij} \|$, for the  embedded Markov chain $J_n$, 
expectations $\| e^{(1)}_{ij}\|$ and second moments $\| e^{(2)}_{ij}\|$  of sojourn times, for the semi-Markov process $J(t)$,   have the 
following forms, 
\begin{equation}\label{partysah}
\left \| 
\begin{array}{cccccc}
0 \ & 0 \  & 0\  & 1 \vspace{1mm} \\
\frac{1}{4} \  & \frac{1}{4} \  & \frac{1}{4} \ & \frac{1}{4} \vspace{1mm}    \\
0 \ & \frac{1}{3} \ & \frac{1}{3} \ & \frac{1}{3}  \vspace{1mm}   \\
0 \  & 0Ê\  & \frac{1}{2} \ & \frac{1}{2}    \\
\end{array}
\right \|, \ 
\left \| 
\begin{array}{cccccc}
0 \ & 0 \  & 0 \ & 1  \vspace{1mm} \\
1 \ & 1 \ & 1 \ & 1  \vspace{1mm} \\
0 \ & 1 \ & 1 \ & 1 \vspace{1mm}  \\
0 \ & 0 \ & 1 \  & 1 \\
\end{array}
\right \| \, {\rm and} \,  
\left \| 
\begin{array}{cccccc}
0 \ & 0 \ & 0 \ & 1  \vspace{1mm} \\
8 \ & 8 \ & 8 \ & 8  \vspace{1mm} \\
0 \ & 6 \ & 6 \ & 6 \vspace{1mm}  \\
0 \ & 0 \ & 2 \ & 2 \\
\end{array}
\right \|. 
\end{equation}

Let us compute first two moments of hitting times $E^{(1)}_{00}, E^{(1)}_{10}$ and $E^{(2)}_{00}, E^{(2)}_{10}$ using the recurrent algorithm  described in Sections 3 -- 5.

Let us first exclude state 3 from the phase space $\XX  = \{ 0, 1, 2, 3 \}$ of the semi-Markov process $J(t)$. 
The corresponding reduced semi-Markov process  $_{\langle 3 \rangle}J(t)$ has the phase space  $_{\langle 3 \rangle}\XX  = \{ 0, 1, 2 \}$. 

The  recurrent formulas (\ref{innmoonan}) and (\ref{dopionmina}) for  transition probabilities of the  embedded Markov chain $_{\langle 3 \rangle}J_n$, 
expectations and second moments of sojourn times  for the semi-Markov process $_{\langle 3 \rangle}J(t)$ have the following forms, respectively, $_{\langle 3 \rangle}p_{ij} = p_{ij} + p_{i3} \frac{p_{3j}}{1 - p_{33}}$,  \ $_{\langle 3 \rangle}e^{(1)}_{ij}  = e^{(1)}_{ij}+  e^{(1)}_{i3} \, _{\langle 3 \rangle}p_{3j}  +  \frac{p_{i3}}{1 - p_{33}} ( e^{(1)}_{3j}+  e^{(1)}_{33} \, _{\langle 3 \rangle}p_{3j})$ and $_{\langle 3 \rangle}e^{(2)}_{ij}  = e^{(2)}_{ij}+   e^{(2)}_{i3} \, _{\langle 3 \rangle}p_{3j} 
+ 2 e^{(1)}_{i3} \,  _{\langle 3 \rangle}e^{(1)}_{3j} 
+  \frac{p_{i3}}{1 - p_{33}} ( e^{(2)}_{3j}+  e^{(2)}_{33} \, _{\langle 3 \rangle}p_{3j}  + 
2 e^{(1)}_{33} \, _{\langle 3 \rangle}e^{(1)}_{3j})$, for $i =  0, 1, 2,  3, \, j = 0, 1, 2$.

The $4 \times 3$ matrices  of transition probabilities $\| _{\langle 3 \rangle}p_{ij} \|$, expectations  $\| _{\langle 3 \rangle}e^{(1)}_{ij}\|$, and second moments $\| _{\langle 3 \rangle}e^{(2)}_{ij}\|$,  computed according the above recurrent formulas,   take the following forms, 
\begin{equation}\label{partysaha}
\left \| 
\begin{array}{ccccc}
0 \ & 0 \ & 1  \vspace{1mm} \\
\frac{1}{4} \  & \frac{1}{4} \  & \frac{1}{2}  \vspace{1mm}    \\
0 \  & \frac{1}{3} \  & \frac{2}{3}      \\
0 \  & 0 \  & 1 
\end{array}
\right \|, \ 
\left \| 
\begin{array}{ccccc}
0 \ & 0 \ & 5   \vspace{1mm} \\
1 \ & 1 \ & 3   \vspace{1mm} \\
0 \ & 1 \ &  \frac{10}{3}    \\
0 \ & 0 \  & 4
\end{array}
\right \| \, {\rm and} \,  
\left \| 
\begin{array}{ccccc}
0 \ & 0 \ & 33   \vspace{1mm} \\
8 \ & 8 \ & 30   \vspace{1mm} \\
0 \ & 6 \ & 28   \\
0 \ & 0 \   & 24
\end{array}
\right \|. 
\end{equation}

Let us now exclude state 2 from the phase space $_{\langle 3 \rangle}\XX = \{ 0, 1, 2 \}$ of the semi-Markov process $_{\langle 3 \rangle}J(t)$.
The corresponding reduced semi-Markov process  $_{\langle 3, 2 \rangle}J(t)$ has the phase space  $_{\langle 3, 2 \rangle}\XX  = \{ 0, 1 \}$.  

The recurrent formulas (\ref{innmoonan}) and (\ref{dopionmina}) for  transition probabilities of the  embedded Markov chain $_{\langle 3, 2 \rangle}J_n$, 
expectations of sojourn times and second moments of sojourn times  for the semi-Markov process $_{\langle 3, 2 \rangle}J(t)$ have the following forms, respectively, $_{\langle 3, 2 \rangle}p_{ij} = \, _{\langle 3 \rangle}p_{ij} + \, _{\langle 3 \rangle}p_{i2} \frac{_{\langle 3 \rangle}p_{2j}}{1 - \, _{\langle 3 \rangle}p_{22}}$,  \ $_{\langle 3, 2 \rangle}e^{(1)}_{ij}  = \, _{\langle 3 \rangle}e^{(1)}_{ij} + \, _{\langle 3 \rangle}e^{(1)}_{i2} \, _{\langle 3, 2 \rangle}p_{2j}  \, +  \, \frac{_{\langle 3 \rangle}p_{i2}}{1 - \, _{\langle 3 \rangle}p_{22}} ( _{\langle 3 \rangle}e^{(1)}_{2j} \, +  \, _{\langle 3 \rangle}e^{(1)}_{22} \, _{\langle 3, 2 \rangle}p_{2j})$ \ and  \
$_{\langle 3, 2 \rangle}e^{(2)}_{ij}  = \, _{\langle 3 \rangle}e^{(2)}_{ij} \, + \,  _{\langle 3 \rangle}e^{(2)}_{i2} \, 
_{\langle 3, 2 \rangle}p_{2j} \, 
+ \, 2 \, _{\langle 3 \rangle}e^{(1)}_{i2} \,  _{\langle 3, 2 \rangle}e^{(1)}_{2j} \,
+   \frac{_{\langle 3 \rangle}p_{i2}}{1 - \, _{\langle 3 \rangle}p_{22}} ( _{\langle 3 \rangle}e^{(2)}_{2j} + \, _{\langle 3 \rangle}e^{(2)}_{22} \, _{\langle 3, 2 \rangle}p_{2j}  + 
2 \, _{\langle 3 \rangle}e^{(1)}_{22} \, _{\langle 3, 2 \rangle}e^{(1)}_{2j})$, for $i =  0, 1, 2, \, j = 0, 1$.

The $3 \times 2$ matrices  of transition probabilities \, $\| _{\langle 3, 2 \rangle}p_{ij} \|$, expectations   $\| _{\langle 3, 2 \rangle}e^{(1)}_{ij}\|$, and second moments $\| _{\langle 3, 2 \rangle}e^{(2)}_{ij}\|$,  computed according the above recurrent formulas,   take the following forms, 
\begin{equation}\label{partysahaku}
 \left \| 
\begin{array}{cccc}
0 \ & 1    \vspace{1mm} \\
\frac{1}{4} \  & \frac{3}{4}     \vspace{1mm}    \\
0 \  & 1        
\end{array}
\right \|, \ 
 \left \| 
\begin{array}{cccc}
0 \ & 18    \vspace{1mm} \\
1 \ & \frac{21}{2}    \vspace{1mm} \\
0 \ &  13    
\end{array}
\right \|  \,   {\rm and} \
\left \| 
\begin{array}{cccc}
0 \ & 525    \vspace{1mm} \\
8 \ & 297    \vspace{1mm} \\
0 \ &  362   
\end{array}
\right \|. 
\end{equation}

Finally, let us  exclude state 1 from the phase space $_{\langle 3, 2 \rangle}\XX = \{ 0, 1 \}$ of the semi-Markov process $_{\langle 3, 2 \rangle}J(t)$.
The corresponding reduced semi-Markov process  $_{\langle 3, 2, 1 \rangle}J(t)$ has the phase space  $_{\langle 3, 2, 1 \rangle}\XX  = \{ 0 \}$. 

The recurrent formulas (\ref{innmoonan}) and (\ref{dopionmina}) for  transition probabilities of the  embedded Markov chain $_{\langle 3, 2, 1 \rangle}J_n$, 
expectations of sojourn times and second moments of sojourn times  for the semi-Markov process $_{\langle 3, 2, 1 \rangle}J(t)$ have the following forms, respectively, $_{\langle 3, 2, 1 \rangle}p_{i0} = \, _{\langle 3, 2 \rangle}p_{i0} + \, _{\langle 3, 2 \rangle}p_{i1} \frac{_{\langle 3, 2 \rangle}p_{10}}{1 - \, _{\langle 3, 2 \rangle}p_{11}}$,   $E^{(1)}_{i0}   = \, _{\langle 3, 2, 1 \rangle}e^{(1)}_{i0}  = \, _{\langle 3, 2 \rangle}e^{(1)}_{i0} + \, _{\langle 3, 2 \rangle}e^{(1)}_{i1} \, _{\langle 3, 2, 1 \rangle}p_{10}  +  \, \frac{_{\langle 3, 2 \rangle}p_{i1}}{1 - \, _{\langle 3, 2 \rangle}p_{1 1}} ( _{\langle 3, 2 \rangle}e^{(1)}_{10}  +  \, _{\langle 3, 2 \rangle}e^{(1)}_{11} \, _{\langle 3, 2, 1 \rangle}p_{10})$ and 
$E^{(2)}_{i0}  = \, _{\langle 3, 2, 1 \rangle}e^{(2)}_{i 0}   = \, _{\langle 3, 2 \rangle}e^{(2)}_{i0} \ + \,  _{\langle 3, 2 \rangle}e^{(2)}_{i1} \, 
_{\langle 3, 2, 1 \rangle}p_{1 0} \ 
+  2 \, _{\langle 3, 2 \rangle}e^{(1)}_{i1} \,  _{\langle 3, 2, 1 \rangle}e^{(1)}_{1 0} \
+  \  \frac{_{\langle 3, 2 \rangle}p_{i1}}{1 - \, _{\langle 3, 2 \rangle}p_{11}} ( _{\langle 3, 2 \rangle}e^{(2)}_{1 0} \ + \, 
_{\langle 3, 2 \rangle}e^{(2)}_{11} \, _{\langle 3, 2, 1 \rangle}p_{10}  + 
2 \, _{\langle 3, 2 \rangle}e^{(1)}_{11} \, _{\langle 3, 2, 1 \rangle}e^{(1)}_{10})$, for $i =  0, 1$.

Here, equalities,  $_{\langle 3, 2, 1 \rangle}e^{(0)}_{i0} = \ _{\langle 3, 2, 1 \rangle}p_{i0} = 1, i = 0, 1$,  should  be taken into account that simplifies the corresponding calculations.

The $2 \times 1$ matrices  of expectations  $\| E^{(1)}_{i 0} \|$, and second moments $\| E^{(2)}_{i 0}\|$  computed according the above recurrent formulas,   take the following forms, 
\begin{equation}\label{partysahamop}
\| E^{(1)}_{i 0} \| = \left \| 
\begin{array}{lllllll}
64    \vspace{1mm} \\
46   
\end{array}
\right \| \ {\rm and} \  
\| E^{(2)}_{i 0} \| = \left \| 
\begin{array}{lllllll}
7265    \vspace{1mm} \\
5084    
\end{array}
\right \|. 
\end{equation}

In conclusion, we would like to note that recurrent algorithms presented in the paper are subjects of effective program realization. These programs let one compute 
power moments for hitting times and accumulated rewards of hitting times for semi-Markov processes with very large numbers of states. We are going to present such programs and results of 
large scale experimental studies in  future publications.


\begin{thebibliography}{99}
{\footnotesize

\bibitem{An1}
Anisimov, V.~V. (2008). Switching Processes in Queueing Models. Applied Stochastic Methods 
Series. ISTE, London and Wiley, Hoboken, NJ, 352 pp. 
\vspace{-2mm}

\bibitem{AZD}
Anisimov, V.~V., Zakusilo, O.~K., Donchenko, V.~S. (1987). Elements
of Queueing and Asymptotical Analysis of Systems. Lybid', Kiev, 246 pp.  
\vspace{-2mm}

\bibitem{ABL}
Barbu, V.,  Boussemart, M.,  Limnios, N.  (2004). Discrete time semi-Markov model for reliability and survival
analysis. Commun. Stat. Theory, Methods, 33, no. 11, 2833--2868.
\vspace{-2mm}

\bibitem{BDB}
Biffi, G., D'amigo, G.,  Di Biase, G.,  Janssen, J., Manca, R., Silvestrov, D. (2008). 
Monte Carlo semi-Markov methods  for credit risk migration and Basel II rules. I, II. J. Numer. 
Appl. Math., 1(96),  I: 28--58, II: 59--86.
\vspace{-2mm}

\bibitem{Chu1}
Chung, K.~ L. (1954). Contributions to the theory
of Markov chains II. Trans. Amer. Math. Soc., 76, 397--419. 
\vspace{-2mm}

\bibitem{Ghu2}
Chung, K.~L.  (1960, 1967). Markov Chains with Stationary
Transition Probabilities. Fundamental Principles of Mathematical
Sciences, 104, Springer, Berlin, x+278 pp.

\vspace{-2mm}

\bibitem{CRST}
Ciardo, G.,  Raymonf, M.~A., Sericola, B., Trivedi,  K.~S. (1990). Performability analysis using semi-Markov reward
processes. IEEE Trans. Comput., 39, no. 10, 1251--1264.
\vspace{-2mm}

\bibitem{Cog}
Cogburn, R.  (1975).  A uniform theory for sums of Markov chain transition probabilities. Ann. Probab., 3, 191--214.
\vspace{-2mm}

\bibitem{Cou}
Courtois, P.~J. 1977). Decomposability. Queueing and Computer System  Applications. ACM Monograph Series, Academic Press, New  York, xiii+201 pp. 
\vspace{-2mm}

\bibitem{DGM}
D'Amico, G., Guillen, M., Manca, R. (2013).  Semi-Markov disability insurance models. Commun. Stat. Theory
Methods, 42, no. 16, 2172--2188.
\vspace{-2mm}

\bibitem{DJM}
D'Amico, G., Janssen, J., Manca, R.  (2005). Homogeneous semi-Markov
reliability models for credit risk management. Decis. Econ. Finance, 28, no. 2, 79--93.
\vspace{-2mm}

\bibitem{DP1}
D'Amico, G., Petroni, F.  (2012).  A semi-Markov model for price returns. Physica A, 391, 4867--4876.
\vspace{-2mm}

\bibitem{DPP1}
D'Amico, G., Petroni, F., Prattico, F. (2013). First and second order semi-Markov chains for wind speed
modeling. Physica A, 392, no. 5, 1194--1201.
\vspace{-2mm}

\bibitem{DPP2}
D'Amico, G., Petroni, F.,  Prattico, F. (2015). Performance analysis of second order semi-Markov chains: an application to wind energy production
Methodol. Comput. Appl. Probab., 17, no. 3, 781--794.
\vspace{-2mm}

\bibitem{GS1}
Gyllenberg, M., Silvestrov, D.  (2008). Quasi-Stationary Phenomena in Nonlinearly Perturbed Stochastic Systems. De
Gruyter Expositions in Mathematics, 44, Walter de Gruyter, Berlin,  ix+579 pp.
\vspace{-2mm}

\bibitem{Hunt2}
Hunter, J.~J.  (2005).  Stationary distributions and mean first passage
times of perturbed Markov chains. Linear Algebra Appl.,  410, 217--243.
\vspace{-2mm}

\bibitem{JM1}
Janssen, J., Manca, R. (2006). Applied Semi-Markov Processes.
Springer, New York, xii+309 pp.  
\vspace{-2mm}

\bibitem{JM2}
Janssen, J., Manca, R. (2007).   Semi-Markov Risk Models for
Finance, Insurance and Reliability.  Springer, New York, xvii+429 pp.
\vspace{-2mm}

\bibitem{KS1}
Kemeny, J.~G., Snell, J.~L. (1961a).  Potentials for
denumerable Markov chains. J. Math. Anal. Appl., 6, 196--260. 
\vspace{-2mm}

\bibitem{KS2}
Kemeny, J.~G., Snell, J.~L.  (1961b).  Finite continuous time Markov chains. Theor. Probab. Appl.,  6, 110--115.
\vspace{-2mm}

\bibitem{KBT}
Korolyuk, V.~S., Brodi, S.~M., Turbin, A.~F. (1974). Semi-Markov processes and their application. Probability Theory. Mathematical 
Statistics. Theoretical Cybernetics, 11, VINTI, Moscow, 1974, 47--97.
\vspace{-2mm}

\bibitem{KoK}
Korolyuk, V.~S., Korolyuk, V.~V. (1999). Stochastic Models of Systems. Mathematics and its Applications, 469, Kluwer, Dordrecht,  xii+185 pp.
\vspace{-2mm}

\bibitem{KoL}
Koroliuk, V.~S., Limnios, N. (2005). Stochastic Systems in
Merging Phase Space. World Scientific, Singapore, xv+331 pp.
\vspace{-2mm}

\bibitem{KT1}
Korolyuk, V.~S., Turbin, A.~F. (1976). Semi-Markov Processes and its Applications. Naukova Dumka, Kiev, 184 pp.
\vspace{-2mm}

\bibitem{KT2}
Korolyuk, V.~S., Turbin, A.~F. (1978). Mathematical Foundations of the State Lumping of Large Systems. Naukova Dumka, Kiev, 218 pp.
(English edition: Mathematics and its Applications,  264, Kluwer, Dordrecht, 1993, x+278 pp.). 
\vspace{-2mm}

\bibitem{Kov}
Kovalenko, I.~N. (1975).  Studies in the Reliability Analysis of Complex Systems. Naukova Dumka, Kiev,  210 pp.
\vspace{-2mm}

\bibitem{KKP}
Kovalenko, I.~N., Kuznetsov, N.~Yu., Pegg, P.~A. (1997).  
Mathematical Theory of Reliability of  Time Dependent Systems with
Practical Applications. Wiley Series in Probability and Statistics,
Wiley, New York, 316 pp.
\vspace{-2mm}

\bibitem{Lam}
Lamperty, J. (1963).  Criteria for stochastic processes II: passage-time moments. J. Math. Anal. Appl., 7,  127--145.
\vspace{-2mm}

\bibitem{LiO1}
Limnios, N., Opri\c{s}an, G. (2001).  Semi-Markov Processes and
Reliability. Statistics for Industry and Technology, 
Birkh\"{a}user, Boston, xii+222 pp.
\vspace{-2mm}

\bibitem{LiO2}
Limnios, N., Opri\c{s}an, G. (2003). An introduction to Semi-Markov processes with application to reliability. In:
Shanbhag, D.N., Rao, C. R. (Eds.) Handbook of Statistics, 21, 515--556. 
\vspace{-2mm}

\bibitem{Num}
Nummelin, E. (1984).  General Irreducible Markov Chains and
Nonnegative Operators. Cambridge Tracts in Mathematics,  83, Cambridge University Press, Cambridge, xi+172 pp.
\vspace{-2mm}

\bibitem{Pa2}
Papadopoulou, A.~A., Tsaklidis, G., McClean, S., Garg, L. (2012). On the moments and the distribution of the cost
of a semi-Markov model for healthcare systems. Methodol. Comput. Appl. Probab., 14, no. 3, 717--737.
\vspace{-2mm}

\bibitem{Pa3}
Papadopoulou, A.~A. (2013). Some results on modeling biological sequences and web navigation with a semi-
Markov chain.  Comm. Stat. Theory, Methods, 42, no. 16, 2153--2171.
\vspace{-2mm}

\bibitem{Pi1}
Pitman, J.~W.  (1974a). An identity for stopping times of a Markov
process. In: Williams, E. J. (Ed.) Studies in Probability and Statistics. Academic Press. Jerusalem, 41--57.

\bibitem{Pi2}
Pitman, J.~W.  (1974b). Uniform rates of convergence for Markov chain transition probabilities. Z. Wahrscheinlichkeitsth., 29, 193--227.
\vspace{-2mm}

\bibitem{Pi3}
Pitman, J.~W.  (1977). Occupation measures for Markov
chains. Adv. Appl. Probab., 9,  69--86.
\vspace{-2mm}

\bibitem{Si1}
Silvestrov, D.~S. (1974). Limit Theorems for Composite Random
Functions. Vysshaya Shkola and Izdatel'stvo Kievskogo Universiteta, Kiev, 318 pp. 
\vspace{-2mm}

\bibitem{Si2}
Silvestrov D.~S.  (1980a). Mean hitting times for semi-Markov processes, and queueing networks. Elektron. Infor. Kybern., 16, 399--415.
\vspace{-2mm}

\bibitem{Si3}
Silvestrov, D.~S. (1980b). Semi-Markov Processes with a Discrete
State Space. Library for an Engineer in Reliability, Sovetskoe
Radio, Moscow,  272 pp.  
\vspace{-2mm}

\bibitem{Si4}
Silvestrov, D.~S.  (1983a). Method of a single probability space in ergodic
theorems for regenerative processes I. Math. Operationsforsch.  Statist., Ser. Optimization, 14, 286--299.
\vspace{-2mm}

\bibitem{Si5}
Silvestrov, D.~S.  (1983b).  Invariance principle for the processes with semi-Markov switch-overs with an arbitrary state space.
In:   Proceedings of the Fourth USSR-Japan Symposium on Probability Theory and Mathematical
Statistics, Tbilisi 1983. Lecture Notes in
Math., 1021, 617--628.
\vspace{-2mm}

\bibitem{Si6}
Silvestrov, D.~S.  (1994).   Coupling for Markov renewal processes and the rate of convergence in ergodic theorems for 
processes with semi-Markov switchings. Acta Applic. Math., 34, 109--124.
\vspace{-2mm}

\bibitem{Si7}
Silvestrov, D~S.  (1996). Recurrence relations for generalised hitting
times for semi-Markov processes. Ann. Appl. Probab., 6, 617--649.
\vspace{-2mm}

\bibitem{SD1}  
Silvestrov, D.~S., Drozdenko, M.~O. (2006). Necessary and sufficient conditions for weak convergence of first-rare-event times for semi-Markov processes.  Theory Stoch. Process., 12(28),  no. 3-4, Part I: 151--186, Part II: 187--202.  
\vspace{-2mm}

\bibitem{SSM2}
Silvestrov, D.,  Manca, R., Silvestrova, E. (2014).  Computational algorithms for moments of accumulated Markov and semi-Markov rewards. Comm. Statist. 
Theory Methods, 43,  no. 7,  1453--1469.
\vspace{-2mm}

\bibitem{SSM}
Silvestrov, D., Silvestrova, E., Manca, R. Stochastically ordered models for credit rating dynamics. J. Numer. Appl.
Math., 1(96), 206--218, (2008). 
\vspace{-2mm}

\bibitem{SS1}  
Silvestrov, D., Silvestrov, S.  (2015). Asymptotic expansions for stationary distributions of perturbed semi-Markov processes. Research Report 2015-9, 
Department of Mathematics, Stockholm University, 75 pp.
\vspace{-2mm}

\bibitem{SMS1}
Stenberg, F., Manca, R., Silvestrov, D. (2006).  Semi-Markov reward models for disability
insurance. Theory Stoch. Proces., 12(28), no. 3-4,  239--254.
\vspace{-2mm}

\bibitem{SMS2}
Stenberg, F., Manca, R., Silvestrov, D. (2007).  An algorithmic approach to discrete time non-homogeneous backward semi-Markov reward process with an application to disability insurance.   Metodol. Comput. Appl. Probab., 9, 497--519.
\vspace{-2mm}

\bibitem{YZ1}
Yin, G.~G., Zhang, Q.  (2005).  Discrete-Time Markov Chains.
Two-Time-Scale Methods and Applications. Stochastic Modelling and Applied Probability, 55, Springer, New York, xix+348 pp.
\vspace{-2mm}

\bibitem{YZ2}
Yin, G.~G., Zhang, Q. (2013). Continuous-Time Markov Chains and Applications. A Two-Time-Scale Approach. Second edition,  Stochastic Modelling and 
Applied Probability, 37, Springer, New York,  xxii+427 pp. (An extended variant of the first (1998)  edition).

}
\end{thebibliography}
\end{document}